\documentclass[12pt,oneside,reqno]{amsart}
\usepackage{mathrsfs}
\usepackage{graphics}
\usepackage{enumerate, amssymb}
\newcommand{\dif}{\mathrm{d}}

\newcommand{\be}{\begin{eqnarray}}
\newcommand{\ee}{\end{eqnarray}}
\newcommand{\ce}{\begin{eqnarray*}}
\newcommand{\de}{\end{eqnarray*}}
\newtheorem{theorem}{Theorem}[section]
\newtheorem{lemma}[theorem]{Lemma}
\newtheorem{remark}[theorem]{Remark}
\newtheorem{definition}[theorem]{Definition}
\newtheorem{proposition}[theorem]{Proposition}
\newtheorem{Example}[theorem]{Example}
\newtheorem{corollary}[theorem]{Corollary}
\def\e{\varepsilon}

\def\a{\alpha}

\def\[{{\Big[}}
\def\]{{\Big]}}
\def\<{{\langle}}
\def\>{{\rangle}}
\def\({{\Big(}}
\def\){{\Big)}}

\def\no{\nonumber}
\def\bt{\begin{theorem}}
\def\et{\end{theorem}}
\def\bl{\begin{lemma}}
\def\el{\end{lemma}}
\def\br{\begin{remark}}
\def\er{\end{remark}}
\def\bx{\begin{Example}}
\def\ex{\end{Example}}
\def\bd{\begin{definition}}
\def\ed{\end{definition}}
\def\bp{\begin{proposition}}
\def\ep{\end{proposition}}
\def\bc{\begin{corollary}}
\def\ec{\end{corollary}}

\def\cB{{\mathcal B}}
\def\cC{{\mathcal C}}

\def\cL{{\mathcal L}}
\def\cM{{\mathcal M}}
\def\cN{{\mathcal N}}

\def\mE{{\mathbb E}}

\def\mN{{\mathbb N}}

\def\mP{{\mathbb P}}

\def\mR{{\mathbb R}}

\def\mU{{\mathbb U}}

\def\geq{\geqslant}
\def\leq{\leqslant}

\begin{document}

\allowdisplaybreaks

\title{Convergence of nonlinear filterings for multiscale systems with correlated L\'evy noises*}

\author{Huijie Qiao}



\thanks{*This work was partly supported by NSF of China (No. 11001051, 11371352) and China Scholarship Council under Grant No. 201906095034.}


\date{}

\dedicatory{School of Mathematics,
Southeast University\\
Nanjing, Jiangsu 211189,  China\\
Department of Mathematics, University of Illinois at
Urbana-Champaign\\
Urbana, IL 61801, USA\\
hjqiaogean@seu.edu.cn}

\begin{abstract}
In the paper, we consider nonlinear filtering problems of multiscale systems in two cases-correlated sensor L\'evy noises and correlated L\'evy noises. First of all, we prove that the slow part of the origin system converges to the homogenized system in the uniform mean square sense. And then based on the convergence result, in the case of correlated sensor L\'evy noises, the nonlinear filtering of the slow part is shown to approximate that of the homogenized system in $L^1$ sense. However, in the case of correlated L\'evy noises, we prove that the nonlinear filtering of the slow part converges weakly to that of the homogenized system.
\end{abstract}

\maketitle \rm

\section{Introduction}

Nowadays, more and more high dimensional and complex mathematical models are used in engineering and science(c.f. \cite{rbn, Dit, kus2, ps1, q00, qzd, yz, zqd}). For example, in some climate models, it is common to simulate the dynamics of the atmosphere and ocean on varying spatial grids with distinct time scale separations.

Simultaneously, controlling, estimating and forecasting these models become more and more interesting(c.f. \cite{rbn, kus2, ps1, qzd, zqd}). However, different time scales make much trouble. Therefore, how to treat these scales is the first important task. A kind of usual methods is to reduce the dimension of these high dimensional mathematical models and study low dimensional ones with the similar dynamical structure. Thus, by estimating the low dimensional ones, we can master the origin high dimensional ones. And nonlinear filtering problems are just right to estimate unobservable and complicated phenomena by observing some simple objects. So, by solving some suitable filtering problems, high dimensional and complex models can be controlled and estimated.

In the paper, we are mainly interested in the nonlinear filtering problems of the following two multiscale systems. For a fixed time $T>0$, given a completed filtered probability space $(\Omega, \mathscr{F}, \{\mathscr{F}_t\}_{t\in[0,T]},\mP)$. Consider the following slow-fast system on $\mR^n\times\mR^m$ and the observation process on $\mR^d$: for $0\leq t\leq T$,
\be\left\{\begin{array}{l}
\dif X^\e_t=b_1(X^\e_t,Z^\e_t)\dif t+\sigma_1(X^\e_t)\dif V_t+\int_{\mU_1}f_1(X^\e_{t-}, u)\tilde{N}_{p_1}(\dif t, \dif u), \\
X^\e_0=x_0,\\
\dif Z^\e_t=\frac{1}{\e}b_2(X^\e_t,Z^\e_t)\dif t+\frac{1}{\sqrt{\e}}\sigma_2(X^\e_t,Z^\e_t)\dif W_t+\int_{\mU_2}f_2(X^\e_{t-},Z^\e_{t-},u)\tilde{N}^{\e}_{p_2}(\dif t, \dif u),\\
Z^\e_0=z_0,\\
\dif Y_t^{\e}=h(X_t^{\e})\dif t+\sigma_3 \dif V_t+\sigma_4 \dif B_t,\\
Y_0^{\e}=0,
\end{array}
\right.
\label{Eq0}
\ee
where $V, W, B$ are $l$-dimensional, $m$-dimensional  and $j$-dimensional standard Brownian motions, respectively, and $p_1, p_2$ are
two stationary Poisson point processes of the class (quasi left-continuous) defined on $(\Omega,
\mathscr{F}, \{\mathscr{F}_t\}_{t\in[0,T]},\mP)$ with values in $\mU$ and the characteristic
measure $\nu_1, \nu_2$, respectively. Here $\nu_1, \nu_2$ are two $\sigma$-finite measures defined
on a measurable space ($\mU,\mathscr{U}$). Fix $\mU_1, \mU_2\in\mathscr{U}$ with $\nu_1(\mU\setminus\mU_1)<\infty$
and $\nu_2(\mU\setminus\mU_2)<\infty$. Let $N_{p_1}((0,t],\dif u)$ be the counting measure of $p_1(t)$, a Poisson random measure
and then $\mE N_{p_1}((0,t],A)=t\nu_1(A)$ for $A\in\mathscr{U}$. Denote
\ce
\tilde{N}_{p_1}((0,t],\dif u):=N_{p_1}((0,t],\dif u)-t\nu_1(\dif u), \qquad\qquad A\in\mathscr{U}|_{\mU_1},
\de
the compensated measure of $N_{p_1}((0,t],\dif u)$. By the same way, we could define $N_{p_2}((0,t],\dif u)$, $\tilde{N}_{p_2}((0,t],\dif u)$.
And $N_{p_2}^{\e}((0,t],\dif u)$ is another Poisson random measure on ($\mU,\mathscr{U}$) such that $\mE N_{p_2}^{\e}((0,t],A)=\frac{1}{\e}t\nu_2(A)$ for $A\in\mathscr{U}$. Moreover, $V_t, W_t, B_t, N_{p_1}, N_{p_2}, N_{p_2}^{\e}$ are mutually independent. The mappings $b_1:\mR^n\times\mR^m\mapsto\mR^n$, $b_2:\mR^n\times\mR^m\mapsto\mR^m$, $\sigma_1:\mR^n\mapsto\mR^{n\times l}$, $\sigma_2:\mR^n\times\mR^m\mapsto\mR^{m\times m}$, $f_1:\mR^n\times\mU_1\mapsto\mR^n$, $f_2:\mR^n\times\mR^m\times\mU_2\mapsto\mR^m$, and $h: \mR^n\mapsto\mR^d$
are all Borel measurable. The matrices $\sigma_3, \sigma_4$ are $d\times l, d\times j$, respectively. The system (\ref{Eq0}) is usually called as a correlated sensor noise model.

We also consider the following slow-fast system on $\mR^n\times\mR^m$ and the observation process on $\mR^d$: for $0\leq t\leq T$, $l=d$,
\be\left\{\begin{array}{l}
\dif \check{X}^\e_t=\check{b}_1(\check{X}^\e_t, \check{Z}^\e_t)\dif t+\check{\sigma}_0(\check{X}^\e_t)\dif B_t+\check{\sigma}_1(\check{X}^\e_t)\dif V_t+\int_{\mU_1}\check{f}_1(\check{X}^\e_{t-}, u)\tilde{N}_{p_1}(\dif t, \dif u), \\
\check{X}^\e_0=\check{x}_0,\\
\dif \check{Z}^\e_t=\frac{1}{\e}\check{b}_2(\check{X}^\e_t,\check{Z}^\e_t)\dif t+\frac{1}{\sqrt{\e}}\check{\sigma}_2(\check{X}^\e_t,\check{Z}^\e_t)\dif W_t+\int_{\mU_2}\check{f}_2(\check{X}^\e_{t-},\check{Z}^\e_{t-},u)\tilde{N}^{\e}_{p_2}(\dif t, \dif u),\\
\check{Z}^\e_0=\check{z}_0,\\
\dif \check{Y}_t^{\e}=\check{h}(\check{X}_t^{\e})\dif t+V_t+\int_0^t\int_{\mU_3}\check{f}_3(s,u)\tilde{N}_{\lambda}(\dif s, \dif u)+\int_0^t\int_{\mU\setminus\mU_3}\check{g}_3(s,u)N_{\lambda}(\dif s, \dif u)\\
\check{Y}_0^{\e}=0
\end{array}
\right.
\label{Eq01}
\ee
where $N_{\lambda}(\dif t,\dif u)$ is a random measure with a predictable compensator $\lambda(t,\check{X}_t^\e,u)\dif t\nu_3(\dif u)$. Here the function $\lambda: [0,T]\times\mR^n\times\mU\rightarrow(0,1)$ is Borel measurable and $\nu_3$ is a $\sigma$-finite measure defined on $\mU$ with $\nu_3(\mU\setminus\mU_3)<\infty$ and $\int_{\mU_3}\|u\|_{\mU}^2\,\nu_3(\dif u)<\infty$ for a fixed $\mU_3\in\mathscr{U}$. Concretely speaking, set 
$$
\tilde{N}_\lambda((0,t], A):=N_\lambda((0,t],A)-\int_0^t\int_A\lambda(s,X_s^\e,u)\dif s\nu_3(\dif u), \quad t\in[0,T], A\in\mathscr{U}|_{\mU_3},
$$ 
and then $\tilde{N}_\lambda((0,t],\dif u)$ is the compensated martingale measure
of $N_{\lambda}((0,t],\dif u)$. Moreover, $V_t, W_t, B_t, N_{p_1}, N_{p_2}, N_{p_2}^{\e}, N_{\lambda}$ are mutually independent. These mappings $\check{b}_1:\mR^n\times\mR^m\mapsto\mR^n$, $\check{b}_2:\mR^n\times\mR^m\mapsto\mR^m$, $\check{\sigma}_0:\mR^n\mapsto\mR^{n\times j}$, $\check{\sigma}_1:\mR^n\mapsto\mR^{n\times d}$, $\check{\sigma}_2:\mR^n\times\mR^m\mapsto\mR^{m\times m}$, $\check{f}_1:\mR^n\times\mU_1\mapsto\mR^n$, $\check{f}_2:\mR^n\times\mR^m\times\mU_2\mapsto\mR^m$, $\check{h}: \mR^n\mapsto\mR^d$, $\check{f}_3:[0,T]\times\mU_3\mapsto\mR^d$ and $\check{g}_3:[0,T]\times(\mU\setminus\mU_3)\mapsto\mR^d$ are all Borel measurable. As usual, the system (\ref{Eq01}) is called as a correlated noise model.

Note that in the systems (\ref{Eq0}) (\ref{Eq01}), the unobservable processes and the observable ones have correlated parts. The type of the multiscale correlated filtering problems usually stems from atmospheric and climatology problems. For example, coupled atmosphere-ocean models provide a multiscale model with fast atmospheric and slow ocean dynamics. In the case of climate prediction, the ocean memory, due to its heat capacity, holds important information. Hence, the improved estimate of the ocean state, which is often the slow component, is of greater interest.

In the paper, we firstly prove that the slow part of a fast-slow system converges to the homogenized system in the uniform mean square sense. And then based on the convergence result, for the system (\ref{Eq0}), the nonlinear filtering of the slow part is shown to approximate that of the homogenized system in $L^1$ sense. But for the system (\ref{Eq01}), we prove that the nonlinear filtering of the slow part converges weakly  to that of the homogenized system.

It is worthwhile to mentioning our method and results. Firstly, for the system (\ref{Eq0}), since the driving processes of the fast-slow system are correlated with that of the observation, we can {\it not} obtain the Zakai equation of the homogenized system (cf. \cite{q3}). Thus, those methods by means of the Zakai equation do {\it not} work (cf. \cite{ImkellerSri, kus, q00}). Therefore, we make use of the exponential martingale to prove the convergence for the filtering of the slow part to that of the homogenized system. 
However, for the system (\ref{Eq01}), we can deduce the Zakai equations of the slow system and the homogenized system, and then show their filtering convergence. Secondly, here we prove the uniform mean square convergence stronger than weak convergence in \cite{rbn, kus2, q00} and convergence in probability in \cite{ps1}. Thirdly, when $f_1=f_2=0$ in (\ref{Eq0}) and $\check{f}_1=\check{f}_2=\check{f}_3=\check{g}_3=0$ in (\ref{Eq01}), two types of multiscale correlated filtering problems have appeared in \cite{rbn} and \cite{lh}, respectively. In \cite{rbn}, when the slow part of the origin system converges to the homogenized system in distribution, Beeson and Namachchivaya only stated that the filtering of the slow part also converges to the filtering of the homogenized system in $L^p$ sense. Regretfully, they didn't prove the result. Here we show the convergence in $L^1$ sense when $f_1\neq 0, f_2\neq 0$. Therefore, our result generalizes the result in a manner. In \cite{lh}, Lucic and Heunis proved that the slow part converges weakly  to the homogenized system, and the filtering of the slow part also converges weakly  to that of the homogenized system. Here, we establish that the slow part converges to the homogenized system in the uniform mean square sense, and the filtering of the slow part converges to that of the homogenized system in $L^1$ sense when $\check{f}_1\neq 0, \check{f}_2\neq 0$. Thus, our result is better. Finally, in \cite{q00}, we considered the nonlinear filtering problem of the system (\ref{Eq01}) with $\check{\sigma}_1=0$. Here, we permit $\check{\sigma}_1\neq 0$. Therefore, our result is more general in some sense.

The paper is arranged as follows. In next section, we consider strong convergence for the fast-slow system. In Section \ref{notset}, we define nonlinear filtering problem and then show that the filtering of the slow part for the system (\ref{Eq0}) converges to that of the  homogenize system. In Section \ref{cornoi}, the filtering of the slow part for the system (\ref{Eq01}) is proved to converge weakly  to that of the  homogenize system. We summarize all the results in Section \ref{con}.

The following convention will be used throughout the paper: $C$ with or without indices will denote different positive constants whose values may change from one place to
another.

\section{Convergence of some processes}\label{conpro}

In the section, we study strong convergence for the fast-slow system (\ref{Eq0}) when $\e\rightarrow0$.

\subsection{A slow-fast system}\label{sfsy}

In the subsection, we introduce slow-fast systems and the existence and uniqueness of their solutions.

Let us consider the system (\ref{Eq0}). First of all, we give out our assumptions and state some related results.

\begin{enumerate}[\bf{Assumption 1.}]
\item
\end{enumerate}
\begin{enumerate}[($\mathbf{H}^1_{b_1, \sigma_1, f_1}$)]
\item For $x_1, x_2\in\mR^n$, $z_1, z_2\in\mR^m$, there exist $L_{b_1}, L_{\sigma_1}, L_{f_1}>0$ such that
\ce
&&|b_1(x_1, z_1)-b_1(x_2, z_2)|^2\leq L_{b_1}(|x_1-x_2|^2+|z_1-z_2|^2),\\
&&\|\sigma_1(x_1)-\sigma_1(x_2)\|^2\leq L_{\sigma_1}|x_1-x_2|^2,\\
&&\int_{\mU_1}|f_1(x_1,u)-f_1(x_2,u)|^2\,\nu_1(\dif u)\leq L_{f_1}|x_1-x_2|^2,
\de
where $|\cdot|$ and $\|\cdot\|$ denote the length of a vector and the Hilbert-Schmidt norm of
a matrix, respectively.
\end{enumerate}

\begin{enumerate}[($\mathbf{H}^2_{b_1, \sigma_1, f_1}$)]
\item For $x\in\mR^n$, $z\in\mR^m$, there exists a $L_{b_1, \sigma_1, f_1}>0$ such that
$$
|b_1(x,z)|^2+\|\sigma_1(x)\|^2+\int_{\mU_1}|f_1(x,u)|^2\nu_1(\dif u)\leq L_{b_1, \sigma_1, f_1}.
$$
\end{enumerate}

\begin{enumerate}[($\mathbf{H}^1_{b_2}$)]
\item (i) $b_2$ is bi-continuous in $(x, z)$,\\
(ii) There exist $L_{b_2}\geq0, \bar{L}_{b_2}>0$ such that
\ce
&&|b_2(x_1, z)-b_2(x_2, z)|\leq L_{b_2}|x_1-x_2|, \qquad\qquad\qquad x_1, x_2\in\mR^n, z\in\mR^m,\\
&&\<z_1-z_2, b_2(x, z_1)-b_2(x, z_2)\>\leq -\bar{L}_{b_2}|z_1-z_2|^2, \qquad x\in\mR^n, z_1, z_2\in\mR^m,
\de
(iii) For $x\in\mR^n$, $z\in\mR^m$, there exists a constant $\bar{\bar{L}}_{b_2}>0$ such that
$$
|b_2(x,z)|\leq \bar{\bar{L}}_{b_2}(1+|x|+|z|).
$$
\end{enumerate}

\begin{enumerate}[($\mathbf{H}^1_{\sigma_2}$)]
\item For $x_1, x_2\in\mR^n$, $z_1, z_2\in\mR^m$, there exists a constant $L_{\sigma_2}>0$ such that
\ce
\|\sigma_2(x_1, z_1)-\sigma_2(x_2, z_2)\|\leq L_{\sigma_2}(|x_1-x_2|+|z_1-z_2|).
\de
\end{enumerate}

\begin{enumerate}[($\mathbf{H}^1_{f_2}$)]
\item There exists a positive function $L(u)$ satisfying 
\ce 
\sup_{u\in\mU_2}L(u)\leq\gamma<1~\mbox{ and } \int_{\mU_2}L(u)^2\,\nu_2(\dif u)<+\infty, 
\de 
such that for any $x_1, x_2\in\mR^n$, $z_1, z_2\in\mR^m$ and $u\in\mU_2$ 
\ce
|f_2(x_1,z_1,u)-f_2(x_2,z_2,u)|\leq L(u)(|x_1-x_2|+|z_1-z_2|),
\de 
and 
\ce 
|f_2(0,0,u)|\leq L(u). 
\de
\end{enumerate}

\medspace

Under {\bf Assumption 1.}, by Theorem 1.2 in \cite{q2}, we know that the system (\ref{Eq0}) has a unique strong solution denoted by $(X^\e_t,Z^\e_t)$. 

\subsection{The fast equation}\label{fas}

In the subsection, we mainly study the second part of the system (\ref{Eq0}). 

First, take any $x\in\mR^n$ and fix it. And consider the following SDE in $\mR^m$:
\ce\left\{\begin{array}{l}
\dif Z^x_t=b_2(x,Z^x_t)\dif t+\sigma_2(x,Z^x_t)\dif W_t+\int_{\mU_2}f_2(x,Z^x_t,u)\tilde{N}_{p_2}(\dif t, \dif u),\\
Z^x_0=z_0, \qquad t\geq0.
\end{array}
\right.
\de
Under the assumption ($\mathbf{H}^1_{b_2}$) ($\mathbf{H}^1_{\sigma_2}$) ($\mathbf{H}^1_{f_2}$), the above equation has a unique solution $Z^x_t$.
In addition, it is a Markov process and its transition probability is denoted by $p(x; z_0,t,A)$ for $t\geq0$ and $A\in\mathscr{B}(\mR^m)$. We assume:
\begin{enumerate}[\bf{Assumption 2.}]
\item
\end{enumerate}
\begin{enumerate}[($\mathbf{H}^2_{\sigma_2}$)]
\item There exists a function $\a_1(x)>0$ such that
\ce
\<\sigma_2(x,z)h,h\>\geq\sqrt{\a_1(x)}|h|^2, \qquad z,h\in\mR^m,
\de
and 
\ce
\|\sigma_{\a_1}(x,z_1)-\sigma_{\a_1}(x,z_2)\|^2\leq L_{\a_1}|z_1-z_2|^2, \qquad z_1, z_2\in\mR^m,
\de 
where $\sigma_{\a_1}(x,z)$ is the unique symmetric nonnegative definite matrix such that $\sigma_{\a_1}(x,z)\sigma_{\a_1}(x,z)
=\sigma_2(x,z)\sigma^T_2(x,z)-\a_1(x)\emph{I}$ for the unit matrix $\emph{I}$.
\end{enumerate}
\begin{enumerate}[($\mathbf{H}^1_{b_2,\sigma_2,f_2}$)] 
\item There exist a $r>2$ and two functions $\a_2(x)>0$, $\a_3(x)\geq0$ such that for all $z\in\mR^m$ 
\ce
2\<z,b_2(x,z)\>+\|\sigma_2(x,z)\|^2+\int_{\mU_2}\big|f_2(x,z,u)\big|^2\nu_2(\dif
u)\leq-\a_2(x)|z|^r+\a_3(x). 
\de
\end{enumerate}
\begin{enumerate}[($\mathbf{H}^2_{b_2,\sigma_2,f_2}$)] 
\item $$
M:=2\bar{L}_{b_2}-L_{b_2}-2L^2_{\sigma_2}-2\int_{\mU_2}L^2(u)\nu_2(\dif u)>0.
$$
\end{enumerate}

Under the assumptions ($\mathbf{H}^1_{b_2}$) ($\mathbf{H}^1_{\sigma_2}$) ($\mathbf{H}^1_{f_2}$) ($\mathbf{H}^2_{\sigma_2}$) ($\mathbf{H}^1_{b_2,\sigma_2,f_2}$), by Theorem 1.3 in  \cite{q1} it holds that there exists a unique invariant probability measure $\bar{p}(x,\cdot)$  for $Z^x_t$ such that 
\be
\|p(x; z_0,t,\cdot)-\bar{p}(x,\cdot)\|_{var}\leq Ce^{-\a t}, \quad t>0,
\label{experg}
\ee
where $\|\cdot\|_{var}$ is the total variance norm and $C, \a>0$ are two constants independent of $z_0, t$.

\subsection{The homogenized equation}\label{ave}

In the subsection, we construct a homogenized equation and study the relationship between the origin equation and the homogenized one.

Next, set
\ce
\bar{b}_1(x):=\int_{\mR^m}b_1(x,z)\bar{p}(x,\dif z), 
\de
and by \cite[Lemma 3.1]{q00}  we know that $\bar{b}_1$ is Lipschitz continuous. So, we construct a SDE on the probability space $(\Omega, \mathscr{F}, \{\mathscr{F}_t\}_{t\in[0,T]}, \mP)$ as follows:
\be\left\{\begin{array}{l}
\dif X^0_t=\bar{b} _1(X^0_t)\dif t+\sigma_1(X^0_t)\dif V_t+\int_{\mU_1}f_1(X^0_{t-}, u)\tilde{N}_{p_1}(\dif t, \dif u),\\
X^0_0=x_0, \qquad\qquad 0\leq t\leq T.
\label{appequ}
\end{array}
\right.
\ee
Based on the assumptions ($\mathbf{H}^1_{b_1, \sigma_1, f_1}$) ($\mathbf{H}^2_{b_1, \sigma_1, f_1}$), it holds that Eq.(\ref{appequ}) has a unique strong solution denoted as $X^0_t$. And then we study the relation between $X^{\e}$ and $X^0$. To do this, we realize a partition of $[0,T]$ into intervals of size $\delta_{\e}>0$, and introduce an auxiliary processes:
\be&&
\dif \hat{Z}^\e_t=\frac{1}{\e}b_2(X^\e_{k\delta_{\e}},\hat{Z}^\e_t)\dif t+\frac{1}{\sqrt{\e}}\sigma_2(X^\e_{k\delta_{\e}},\hat{Z}^\e_t)\dif W_t+\int_{\mU_2}f_2(X^\e_{k\delta_{\e}},\hat{Z}^\e_{t-},u)\tilde{N}^{\e}_{p_2}(\dif t, \dif u),\no\\
&&\qquad\qquad\qquad\qquad\qquad\qquad\qquad    t\in[k\delta_{\e},(k+1)\delta_{\e}),\no\\
&& \hat{Z}^\e_{k\delta_{\e}}=Z^\e_{k\delta_{\e}},
\label{auxpro}
\ee
for $k=0,\cdots, [\frac{T}{\delta_{\e}}]$, where $[\frac{T}{\delta_{\e}}]$ denotes the integer part of $\frac{T}{\delta_{\e}}$. Moreover, we mention the fact that $[\frac{t}{\delta_{\e}}]=k$ for $ t\in[k\delta_{\e},(k+1)\delta_{\e})$. The following lemma gives the relationship between $Z^\e$ and $\hat{Z}^\e$.

\bl
Under {\bf Assumption 1.-2.}, it holds that
\be
\sup\limits_{0\leq s\leq T}\mE|Z^\e_s-\hat{Z}^\e_s|^2\leq \frac{L_{b_2}+2L^2_{\sigma_2}+2\int_{\mU_2}L^2(u)\nu_2(\dif u)}{\e}3(\delta_{\e}+1)L_{b_1, \sigma_1, f_1}\delta^2_{\e}.
\label{zhatz}
\ee
\el
\begin{proof}
By the equations (\ref{Eq0})(\ref{auxpro}), it holds that for $s\in[k\delta_{\e},(k+1)\delta_{\e})$
\ce
Z^\e_s-\hat{Z}^\e_s&=&\frac{1}{\e}\int_{k\delta_{\e}}^s \(b_2(X^\e_r,Z^\e_r)-b_2(X^\e_{k\delta_{\e}},\hat{Z}^\e_r)\)\dif r\\
&&+\frac{1}{\sqrt{\e}}\int_{k\delta_{\e}}^s\(\sigma_2(X^\e_r,Z^\e_r)-\sigma_2(X^\e_{k\delta_{\e}},\hat{Z}^\e_r)\)\dif W_r\\
&&+\int_{k\delta_{\e}}^s\int_{\mU_2}\(f_2(X^\e_r,Z^\e_{r-},u)-f_2(X^\e_{k\delta_{\e}},\hat{Z}^\e_{r-},u)\)\tilde{N}^{\e}_{p_2}(\dif r, \dif u).
\de
Applying the It\^o formula to $Z^\e_s-\hat{Z}^\e_s$ and taking the expectation on two sides, we have that
\ce
\mE|Z^\e_s-\hat{Z}^\e_s|^2&=& \frac{2}{\e}\mE\int_{k\delta_{\e}}^s \<Z^\e_r-\hat{Z}^\e_r, b_2(X^\e_r,Z^\e_r)-b_2(X^\e_{k\delta_{\e}},\hat{Z}^\e_r)\>\dif r\no\\
&&+\frac{1}{\e}\mE\int_{k\delta_{\e}}^s\|\sigma_2(X^\e_r,Z^\e_r)-\sigma_2(X^\e_{k\delta_{\e}},\hat{Z}^\e_r)\|^2\dif r\no\\
&&+\frac{1}{\e}\mE\int_{k\delta_{\e}}^s\int_{\mU_2}|f_2(X^\e_r,Z^\e_r,u)-f_2(X^\e_{k\delta_{\e}},\hat{Z}^\e_r,u)|^2\nu_2(\dif u)\dif r\no\\
&\leq& \frac{2}{\e}\mE\int_{k\delta_{\e}}^s \<Z^\e_r-\hat{Z}^\e_r, b_2(X^\e_r,Z^\e_r)-b_2(X^\e_{k\delta_{\e}},Z^\e_r)\>\dif r\no\\
&&+\frac{2}{\e}\mE\int_{k\delta_{\e}}^s \<Z^\e_r-\hat{Z}^\e_r, b_2(X^\e_{k\delta_{\e}},Z^\e_r)-b_2(X^\e_{k\delta_{\e}},\hat{Z}^\e_r)\>\dif r\no\\
&&+\frac{2L^2_{\sigma_2}}{\e}\mE\int_{k\delta_{\e}}^s\(|X^\e_r-X^\e_{k\delta_{\e}}|^2+|Z^\e_r-\hat{Z}^\e_r|^2\)\dif r\no\\
&&+\frac{2}{\e}\int_{\mU_2}L^2(u)\nu_2(\dif u)\mE\int_{k\delta_{\e}}^s\(|X^\e_r-X^\e_{k\delta_{\e}}|^2+|Z^\e_r-\hat{Z}^\e_r|^2\)\dif r\no\\
&\leq& \frac{2}{\e}\mE\int_{k\delta_{\e}}^s |Z^\e_r-\hat{Z}^\e_r| |b_2(X^\e_r,Z^\e_r)-b_2(X^\e_{k\delta_{\e}},Z^\e_r)|\dif r\no\\
&&-\frac{2\bar{L}_{b_2}}{\e}\mE\int_{k\delta_{\e}}^s |Z^\e_r-\hat{Z}^\e_r|^2\dif r\no\\
&&+\frac{2L^2_{\sigma_2}}{\e}\mE\int_{k\delta_{\e}}^s\(|X^\e_r-X^\e_{k\delta_{\e}}|^2+|Z^\e_r-\hat{Z}^\e_r|^2\)\dif r\no\\
&&+\frac{2}{\e}\int_{\mU_2}L^2(u)\nu_2(\dif u)\mE\int_{k\delta_{\e}}^s\(|X^\e_r-X^\e_{k\delta_{\e}}|^2+|Z^\e_r-\hat{Z}^\e_r|^2\)\dif r\no\\
&\leq& \frac{L_{b_2}}{\e}\mE\int_{k\delta_{\e}}^s (|Z^\e_r-\hat{Z}^\e_r|^2+|X^\e_r-X^\e_{k\delta_{\e}}|^2)\dif r\no\\
&&-\frac{2\bar{L}_{b_2}}{\e}\mE\int_{k\delta_{\e}}^s |Z^\e_r-\hat{Z}^\e_r|^2\dif r\no\\
&&+\frac{2L^2_{\sigma_2}}{\e}\mE\int_{k\delta_{\e}}^s\(|X^\e_r-X^\e_{k\delta_{\e}}|^2+|Z^\e_r-\hat{Z}^\e_r|^2\)\dif r\no\\
&&+\frac{2}{\e}\int_{\mU_2}L^2(u)\nu_2(\dif u)\mE\int_{k\delta_{\e}}^s\(|X^\e_r-X^\e_{k\delta_{\e}}|^2+|Z^\e_r-\hat{Z}^\e_r|^2\)\dif r,
\de
where ($\mathbf{H}^1_{b_2}$) ($\mathbf{H}^1_{\sigma_2}$) ($\mathbf{H}^1_{f_2}$) are used. And then
\ce
\mE|Z^\e_s-\hat{Z}^\e_s|^2+\frac{M}{\e}\mE\int_{k\delta_{\e}}^s|Z^\e_r-\hat{Z}^\e_r|^2\dif r\leq \frac{L_{b_2}+2L^2_{\sigma_2}+2\int_{\mU_2}L^2(u)\nu_2(\dif u)}{\e}\mE\int_{k\delta_{\e}}^s|X^\e_r-X^\e_{k\delta_{\e}}|^2\dif r.
\de
Thus, by ($\mathbf{H}^2_{b_2,\sigma_2,f_2}$) it holds that
\be
\mE|Z^\e_s-\hat{Z}^\e_s|^2&\leq& \frac{L_{b_2}+2L^2_{\sigma_2}+2\int_{\mU_2}L^2(u)\nu_2(\dif u)}{\e}\mE\int_{k\delta_{\e}}^s|X^\e_r-X^\e_{k\delta_{\e}}|^2\dif r.
\label{zhates}
\ee

To obtain (\ref{zhatz}), we only need to estimate $\mE|X^\e_r-X^\e_{k\delta_{\e}}|^2$ for $r\in[k\delta_{\e},(k+1)\delta_{\e})$. Note that 
\ce
X^\e_r-X^\e_{k\delta_{\e}}=\int_{k\delta_{\e}}^r b_1(X^\e_v,Z^\e_v)\dif v+\int_{k\delta_{\e}}^r \sigma_1(X^\e_v)\dif V_v+\int_{k\delta_{\e}}^r\int_{\mU_1}f_1(X^\e_{v-}, u)\tilde{N}_{p_1}(\dif v, \dif u).
\de
So, by the H\"older inequality and ($\mathbf{H}^2_{b_1, \sigma_1, f_1}$) we obtain that
\be
\mE|X^\e_r-X^\e_{k\delta_{\e}}|^2&\leq& 3\mE\left|\int_{k\delta_{\e}}^rb_1(X^\e_v,Z^\e_v)\dif v\right|^2+3\mE\left|\int_{k\delta_{\e}}^r \sigma_1(X^\e_v)\dif V_v\right|^2\no\\
&&+3\mE\left|\int_{k\delta_{\e}}^r\int_{\mU_1}f_1(X^\e_{v-}, u)\tilde{N}_{p_1}(\dif v, \dif u)\right|^2\no\\
&\leq& 3(r-k\delta_{\e})\mE\int_{k\delta_{\e}}^r\left|b_1(X^\e_v,Z^\e_v)\right|^2\dif v+3\mE\int_{k\delta_{\e}}^r \|\sigma_1(X^\e_v)\|^2\dif v\no\\
&&+3\mE\int_{k\delta_{\e}}^r\int_{\mU_1}\left|f_1(X^\e_{v-}, u)\right|^2\nu_1(\dif u)\dif v\no\\
&\leq&3(\delta_{\e}+1)L_{b_1, \sigma_1, f_1}\delta_{\e}.
\label{xdex}
\ee
By inserting (\ref{xdex}) in (\ref{zhates}), it holds that
\ce
\mE|Z^\e_s-\hat{Z}^\e_s|^2&\leq&\frac{L_{b_2}+2L^2_{\sigma_2}+2\int_{\mU_2}L^2(u)\nu_2(\dif u)}{\e}3(\delta_{\e}+1)L_{b_1, \sigma_1, f_1}\delta^2_{\e}.
\de
This is just right (\ref{zhatz}). Thus, the proof is complete.
\end{proof}

Next, we apply (\ref{zhatz}) to estimate $|X^\e_t-X^0_t|$. The main result in the section is the following theorem.

\bt\label{xzerx}
Suppose that {\bf Assumption 1.-2.} hold. Then there exists a constant $C\geq 0$ independent of $\e, \delta_{\e}$ such that 
\be
\mE\(\sup\limits_{0\leq t\leq T}|X^\e_t-X^0_t|^2\)\leq \(C\frac{\e}{\delta_{\e}}+C(\delta_{\e}+1)\delta_{\e}+C(\delta_{\e}+1)\frac{\delta^2_{\e}}{\e}\)e^{CT}.
\label{xzerxe}
\ee
\et
\begin{proof}
By the equations (\ref{Eq0})(\ref{appequ}), we know that
\ce
X^\e_t-X^0_t&=&\int_0^t\left(b_1(X^\e_s,Z^\e_s)-\bar{b}_1(X^0_s)\right)\dif s+\int_0^t\left(\sigma_1(X^\e_s)-\sigma_1(X^0_s)\right)\dif V_s\\
&&+\int_0^t\int_{\mU_1}\left(f_1(X^\e_{s-}, u)-f_1(X^0_{s-}, u)\right)\tilde{N}_{p_1}(\dif s, \dif u), \qquad t\in[0,T].
\de
And then by the Burkholder-Davis-Gundy inequality and the H\"older inequality, it holds that
\be
\mE\(\sup\limits_{0\leq t\leq T}|X^\e_t-X^0_t|^2\)&\leq& 3\mE\(\sup\limits_{0\leq t\leq T}\left|\int_0^t\left(b_1(X^\e_s,Z^\e_s)-\bar{b}_1(X^0_s)\right)\dif s\right|^2\)\no\\
&&+3\mE\(\sup\limits_{0\leq t\leq T}\left|\int_0^t\left(\sigma_1(X^\e_s)-\sigma_1(X^0_s)\right)\dif V_s\right|^2\)\no\\
&&+3\mE\(\sup\limits_{0\leq t\leq T}\left|\int_0^t\int_{\mU_1}\left(f_1(X^\e_{s-}, u)-f_1(X^0_{s-}, u)\right)\tilde{N}_{p_1}(\dif s, \dif u)\right|^2\)\no\\
&\leq&12\mE\(\sup\limits_{0\leq t\leq T}\left|\int_0^t\left(b_1(X^\e_s,Z^\e_s)-b_1(X^\e_{k\delta_{\e}},\hat{Z}^\e_s)\right)\dif s\right|^2\)\no\\
&&+12\mE\(\sup\limits_{0\leq t\leq T}\left|\int_0^t\left(b_1(X^\e_{k\delta_{\e}},\hat{Z}^\e_s)-\bar{b}_1(X^{\e}_{k\delta_{\e}})\right)\dif s\right|^2\)\no\\
&&+12\mE\(\sup\limits_{0\leq t\leq T}\left|\int_0^t\left(\bar{b}_1(X^{\e}_{k\delta_{\e}})-\bar{b}_1(X^{\e}_s)\right)\dif s\right|^2\)\no\\
&&+12\mE\(\sup\limits_{0\leq t\leq T}\left|\int_0^t\left(\bar{b}_1(X^{\e}_s)-\bar{b}_1(X^0_s)\right)\dif s\right|^2\)\no\\
&&+12\mE\int_0^T\left\|\sigma_1(X^\e_s)-\sigma_1(X^0_s)\right\|^2\dif s\no\\
&&+12\mE\int_0^T\int_{\mU_1}\left|f_1(X^\e_{s-}, u)-f_1(X^0_{s-}, u)\right|^2\nu_1(\dif u)\dif s\no\\
&\leq&12 TL_{b_1}\int_0^T\(\mE|X^{\e}_s-X^{\e}_{k\delta_{\e}}|^2+\mE|Z^\e_s-\hat{Z}^\e_s|^2\)\dif s\no\\
&&+12\mE\(\sup\limits_{0\leq t\leq T}\left|\int_0^t\left(b_1(X^\e_{k\delta_{\e}},\hat{Z}^\e_s)-\bar{b}_1(X^\e_{k\delta_{\e}})\right)\dif s\right|^2\)\no\\
&&+12TC\int_0^T\mE|X^{\e}_{k\delta_{\e}}-X^{\e}_s|^2\dif s\no\\
&&+\(12TC+12L_{\sigma_1}+12L_{f_1}\)\int_0^T\mE|X^\e_s-X^0_s|^2\dif s\no\\
&\leq&12\mE\(\sup\limits_{0\leq t\leq T}\left|\int_0^t\left(b_1(X^\e_{k\delta_{\e}},\hat{Z}^\e_s)-\bar{b}_1(X^\e_{k\delta_{\e}})\right)\dif s\right|^2\)\no\\
&&+\(12 TL_{b_1}+12TC\)\int_0^T\mE|X^{\e}_{k\delta_{\e}}-X^{\e}_s|^2\dif s\no\\
&&+12 TL_{b_1}\int_0^T\mE|Z^{\e}_s-\hat{Z}^{\e}_s|^2\dif s\no\\
&&+\(12TC+12L_{\sigma_1}+12L_{f_1}\)\int_0^T\mE\(\sup\limits_{0\leq r\leq s}|X^{\e}_r-X^0_r|^2\)\dif s\no\\
&=:&I_1+I_2+I_3+I_4,
\label{xhatzerest0}
\ee
where ($\mathbf{H}^1_{b_1, \sigma_1, f_1}$) is used in the third inequality. 

Next, we estimate $I_1$. Note that
\be
I_1&=&12\mE\left(\sup\limits_{0\leq i\leq [T/\delta_{\e}]-1}\left|\sum_{k=0}^i\int_{k\delta_{\e}}^{(k+1)\delta_{\e}}\left(b_1(X^\e_{k\delta_{\e}},\hat{Z}^\e_s)-\bar{b}_1(X^\e_{k\delta_{\e}})\right)\dif s\right|^2\right)\no\\
&\leq&12\mE\left(\sup\limits_{0\leq i\leq [T/\delta_{\e}]-1}(i+1)\sum_{k=0}^i\left|\int_{k\delta_{\e}}^{(k+1)\delta_{\e}}\left(b_1(X^\e_{k\delta_{\e}},\hat{Z}^\e_s)-\bar{b}_1(X^\e_{k\delta_{\e}})\right)\dif s\right|^2\right)\no\\
&\leq&12[T/\delta_{\e}]\sum_{k=0}^{[T/\delta_{\e}]-1}\mE\left|\int_{k\delta_{\e}}^{(k+1)\delta_{\e}}\left(b_1(X^\e_{k\delta_{\e}},\hat{Z}^\e_s)-\bar{b}_1(X^\e_{k\delta_{\e}})\right)\dif s\right|^2\no\\
&\leq&12[T/\delta_{\e}]^2\sup\limits_{0\leq k\leq [T/\delta_{\e}]-1}\mE\left|\int_{k\delta_{\e}}^{(k+1)\delta_{\e}}\left(b_1(X^\e_{k\delta_{\e}},\hat{Z}^\e_s)-\bar{b}_1(X^\e_{k\delta_{\e}})\right)\dif s\right|^2\no\\
&\leq&12\(\frac{T}{\delta_{\e}}\)^2\sup\limits_{0\leq k\leq [T/\delta_{\e}]-1}\mE\left|\int_0^{\delta_{\e}}\left(b_1(X^\e_{k\delta_{\e}},\hat{Z}^\e_{k\delta_{\e}+s})-\bar{b}_1(X^\e_{k\delta_{\e}})\right)\dif s\right|^2.
\label{i1est}
\ee
So, we only need to analysis $\mE\left|\int_0^{\delta_{\e}}\left(b_1(X^\e_{k\delta_{\e}},\hat{Z}^\e_{k\delta_{\e}+s})-\bar{b}_1(X^\e_{k\delta_{\e}})\right)\dif s\right|^2$ for $k=0,\cdots, [T/\delta_{\e}]-1$. Fix $k$ and set
\ce\left\{\begin{array}{l}
\dif \check{Z}^\e_t=b_2(X^\e_{k\delta_{\e}},\check{Z}^\e_t)\dif t+\sigma_2(X^\e_{k\delta_{\e}},\check{Z}^\e_t)\dif \check{W}_t+\int_{\mU_2}f_2(X^\e_{k\delta_{\e}},\check{Z}^\e_{t-},u)\tilde{N}_{\check{p}_2}(\dif t, \dif u), t\in[0, \delta_{\e}/\e),\\
\check{Z}^\e_0=Z^\e_{k\delta_{\e}},
\end{array}
\right.
\de
where $\check{W}$, $W$, $\check{p}_2$ and $p_2$ are mutually independent, and $\check{W}$, $W$ and $\check{p}_2$, $p_2$ have the same distributions, respectively. And by the scaling property of Brownian motions and Poission random measures, it holds that $\hat{Z}^\e_{k\delta_{\e}+t}$ and $\check{Z}^\e_{t/\e}$ have the same distribution. Thus we have 
\be
\mE\left|\int_0^{\delta_{\e}}\left(b_1(X^\e_{k\delta_{\e}},\hat{Z}^\e_{k\delta_{\e}+s})-\bar{b}_1(X^\e_{k\delta_{\e}})\right)\dif s\right|^2&=&\mE\left|\int_0^{\delta_{\e}}\left(b_1(X^\e_{k\delta_{\e}},\check{Z}^\e_{s/\e})-\bar{b}_1(X^\e_{k\delta_{\e}})\right)\dif s\right|^2\no\\
&=&{\e}^2\mE\left|\int_0^{{\delta_{\e}}/\e}\left(b_1(X^\e_{k\delta_{\e}},\check{Z}^\e_s)-\bar{b}_1(X^\e_{k\delta_{\e}})\right)\dif s\right|^2\no\\
&=&{\e}^2\mE\int_0^{{\delta_{\e}}/\e}\int_0^{{\delta_{\e}}/\e}\left(b_1(X^\e_{k\delta_{\e}},\check{Z}^\e_r)-\bar{b}_1(X^\e_{k\delta_{\e}})\right)\no\\
&&\qquad\qquad \left(b_1(X^\e_{k\delta_{\e}},\check{Z}^\e_s)-\bar{b}_1(X^\e_{k\delta_{\e}})\right)\dif s\dif r\no\\
&=&2{\e}^2\int_0^{{\delta_{\e}}/\e}\int_r^{{\delta_{\e}}/\e}\mE\left(b_1(X^\e_{k\delta_{\e}},\check{Z}^\e_r)-\bar{b}_1(X^\e_{k\delta_{\e}})\right)\no\\
&&\qquad \left(b_1(X^\e_{k\delta_{\e}},\check{Z}^\e_s)-\bar{b}_1(X^\e_{k\delta_{\e}})\right)\dif s\dif r.
\label{maxest}
\ee
And then we investigate the integrand of the above integration. By the H\"older inequality it holds that
\be
&&\mE\left(b_1(X^\e_{k\delta_{\e}},\check{Z}^\e_r)-\bar{b}_1(X^\e_{k\delta_{\e}})\right) \left(b_1(X^\e_{k\delta_{\e}},\check{Z}^\e_s)-\bar{b}_1(X^\e_{k\delta_{\e}})\right)\no\\
&=&\mE\left[\left(b_1(X^\e_{k\delta_{\e}},\check{Z}^\e_r)-\bar{b}_1(X^\e_{k\delta_{\e}})\right) \mE\left[\left(b_1(X^\e_{k\delta_{\e}},\check{Z}^\e_s)-\bar{b}_1(X^\e_{k\delta_{\e}})\right)|\mathscr{F}^{\check{Z}^\e}_r\right]\right]\no\\
&=&\mE\left[\left(b_1(X^\e_{k\delta_{\e}},\check{Z}^\e_r)-\bar{b}_1(X^\e_{k\delta_{\e}})\right) \mE^{\check{Z}^\e_r}\left(b_1(X^\e_{k\delta_{\e}},\check{Z}^\e_{s-r})-\bar{b}_1(X^\e_{k\delta_{\e}})\right)\right]\no\\
&\leq&\left(\mE\left(b_1(X^\e_{k\delta_{\e}},\check{Z}^\e_r)-\bar{b}_1(X^\e_{k\delta_{\e}})\right)^2\right)^{1/2}\left(\mE\(\mE^{\check{Z}^\e_r}\left(b_1(X^\e_{k\delta_{\e}},\check{Z}^\e_{s-r})-\bar{b}_1(X^\e_{k\delta_{\e}})\right)\)^2\right)^{1/2}\no\\
&\leq&Ce^{-\a(s-r)},
\label{intest}
\ee
where the last inequality is based on ($\mathbf{H}^2_{b_1, \sigma_1, f_1}$) and (\ref{experg}), $\mathscr{F}_r^{\check{Z}^\e} \triangleq\sigma(\check{Z}_v^\e:
 0\leq v \leq r) \vee \cN$ and $\cN$ is the collection of all $\mP$-measure zero sets. Inserting (\ref{intest}) in (\ref{maxest}), we furthermore obtain that 
\be
\mE\left|\int_0^{\delta_{\e}}\left(b_1(X^\e_{k\delta_{\e}},\hat{Z}^\e_{k\delta_{\e}+s})-\bar{b}_1(X^\e_{k\delta_{\e}})\right)\dif s\right|^2&\leq&2{\e}^2\int_0^{{\delta_{\e}}/\e}\int_r^{{\delta_{\e}}/\e}Ce^{-\a(s-r)}\dif s\dif r\no\\
&\leq& C{\e}^2\frac{\delta_{\e}}{\e}.
\label{est}
\ee
By combining (\ref{est}) with (\ref{i1est}), it holds that 
\be
I_1\leq \frac{C}{\delta_{\e}/\e}.
\label{i1este}
\ee

Finally, applying (\ref{i1este}) (\ref{xdex}) (\ref{zhatz}) to (\ref{xhatzerest0}), we have that
\ce
\mE\(\sup\limits_{0\leq t\leq T}|X^\e_t-X^0_t|^2\)&\leq&\frac{C}{\delta_{\e}/\e}+C(\delta_{\e}+1)\delta_{\e}+C(\delta_{\e}+1)\delta^2_{\e}/\e+C\int_0^T\mE\(\sup\limits_{0\leq r\leq s}|X^{\e}_r-X^0_r|^2\)\dif s.
\de
The Gronwall inequality admits us to obtain that
\ce
\mE\(\sup\limits_{0\leq t\leq T}|X^\e_t-X^0_t|^2\)&\leq&\(C\frac{\e}{\delta_{\e}}+C(\delta_{\e}+1)\delta_{\e}+C\frac{(\delta_{\e}+1)\delta^2_{\e}}{\e}\)e^{CT}.
\de
The proof is complete.
\end{proof}

\br\label{del}
Based on Theorem \ref{xzerx}, it holds that $X^\e_t$ converges to $X^0_t$ in the mean square sense if $\frac{\e}{\delta_{\e}}\rightarrow0$ and $\frac{\delta^2_{\e}}{\e}\rightarrow0$ as $\e\rightarrow0$. For example, we take $\delta_{\e}=\e^{2/3}$, and have that $\frac{\e}{\delta_{\e}}=\e^{1/3}\rightarrow0$, and $\frac{\delta^2_{\e}}{\e}=\e^{1/3}\rightarrow0$ when $\e\rightarrow0$. 
\er

\section{Convergence of nonlinear filterings with correlated sensor noises}\label{notset}

In the section, we introduce the nonlinear filtering problems for $X_t^{\e}$ and $X_t^{0}$ and their relationship.

\subsection{Nonlinear filtering problems with the system (\ref{Eq0})} 

In the subsection, we introduce nonlinear filtering problems of $X_t^{\e}$ and $X_t^{0}$.
 
For
\ce
Y_t^{\e}&=&\int_0^th(X_s^{\e})\dif s+\sigma_3 V_t+\sigma_4 B_t,
\de
we make the following hypotheses:
\begin{enumerate}[\bf{Assumption 3.}]
\item
\end{enumerate}
\begin{enumerate}[($\mathbf{H}_{h}$)] 
\item $h$ is bounded.
\end{enumerate}
\begin{enumerate}[($\mathbf{H}_{\sigma_3,\sigma_4}$)] 
\item 
$\sigma_3\sigma^{\prime}_3+\sigma_4\sigma^{\prime}_4=I,$ where $\sigma^{\prime}_3$ stands for the transpose of the matrix $\sigma_3$ and $I$ is the $d$ order unit matrix.
\end{enumerate}

By ($\mathbf{H}_{\sigma_3,\sigma_4}$), we know that $U_t:=\sigma_3 V_t+\sigma_4 B_t$ is a $d$ dimensional Brownian motion. Denote
\ce
(\gamma^\e_t)^{-1}:=\exp\bigg\{-\int_0^t h^i(X^\e_s)\dif U^i_s-\frac{1}{2}\int_0^t
\left|h(X^\e_s)\right|^2\dif s\bigg\}.
\de
Here and hereafter, we use the convention that repeated indices imply summation. And then by ($\mathbf{H}_{h}$) we know that $(\gamma^\e_t)^{-1}$ is an exponential martingale. Define a measure $\mP^\e$ via
$$
\frac{\dif \mP^\e}{\dif \mP}=(\gamma^\e_T)^{-1}.
$$
By the Girsanov theorem for Brownian motions, one can obtain that 
\be\label{tilw}
Y^\e_t=U_t+\int_0^t h(X^\e_s)\dif s
\ee
is a $\mathscr{F}_t$-Brownian motion under the probability measure $\mP^\e$. 

Next, we rewrite $\gamma^\e_t$ as
\ce
\gamma^\e_t&=&\exp\bigg\{\int_0^th^i(X_s^{\e})\dif Y^{\e,i}_s-\frac{1}{2}\int_0^t
\left|h(X_s^{\e})\right|^2\dif s\bigg\}.
\de
Define
\ce
&&\rho^{\e}_t(\psi):=\mE^{\mP^\e}[\psi(X^{\e}_t)\gamma^\e_t|\mathscr{F}_t^{Y^{\e}}], \\
&&\pi^{\e}_t(\psi):=\mE[\psi(X^{\e}_t)|\mathscr{F}_t^{Y^{\e}}], \qquad \psi\in\cB(\mR^n),
\de
where $\mE^{\mP^\e}$ denotes the expectation under the measure $\mP^\e$, $\mathscr{F}_t^{Y^\e} \triangleq\sigma(Y_s^\e:
 0\leq s \leq t) \vee \cN$, $\cN$ is the collection of all $\mP$-measure zero sets  and $\cB(\mR^n)$ denotes the collection of all bounded and Borel measurable functions on $\mR^n$. $ \rho_t^\e$ and $\pi^{\e}_t$ are called the nonnormalized filtering and the normalized filtering of $X_t^\e$ with respect to $\mathscr{F}_t^{Y^{\e}}$, respectively. And then by the Kallianpur-Striebel formula it holds that
\ce
\pi^{\e}_t(\psi)=\frac{\rho^{\e}_t(\psi)}{\rho^{\e}_t(1)}.
\de

Set
\ce
\gamma^0_t&:=&\exp\bigg\{\int_0^th^i(X_s^0)\dif Y^{\e,i}_s-\frac{1}{2}\int_0^t
\left|h(X_s^0)\right|^2\dif s\bigg\},
\de
and furthermore
\ce
\rho^0_t(\psi)&:=&\mE^{\mP^\e}[\psi(X^0_t)\gamma^0_t|\mathscr{F}_t^{Y^{\e}}],\\
\pi^0_t(\psi)&:=&\frac{\rho^0_t(\psi)}{\rho^0_t(1)}.
\de
And then we will prove that $\pi^0$ could be understood as the nonlinear filtering problem for $X_t^0$ with respect to $\mathscr{F}_t^{Y^{\e}}$.

\subsection{The relation between $\pi^{\e}_t$ and $\pi^0_t$}

In the subsection we will show that $\pi^{\e}_t$ converges to $\pi^0_t$ as $\e\rightarrow0$ in a suitable sense. Let us start with a key lemma.

\bl\label{es}
Under $({\bf H}_h)$, there exists a constant $C>0$ such that
$$
\mE\left|\rho^0_t(1)\right|^{-p}\leq\exp\left\{(2p^2+p+1)CT/2\right\}, \quad t\in[0,T], \quad
p>1.
$$
\el
\begin{proof}
Although the proof is similar to Lemma 4.1 in \cite{qzd}, we prove it to the readers' convenience. For $\mE\left|\rho^0_t(1)\right|^{-p}$, we compute
$$
\mE\left|\rho^0_t(1)\right|^{-p}=\mE^\e\left|\rho^0_t(1)\right|^{-p}\gamma^\e_T
\leq(\mE^\e\left|\rho^0_t(1)\right|^{-2p})^{1/2}(\mE^\e(\gamma^\e_T)^2)^{1/2},
$$
where the last inequality is based on the H\"older inequality. For $\mE^\e\left|\rho^0_t(1)\right|^{-2p}$, note that
$\rho^0_t(1)=\mE^\e[\gamma^0_t|\mathscr{F}_t^{Y^\e}]$. And then it follows from the
Jensen inequality that
$$
\mE^\e\left|\rho^0_t(1)\right|^{-2p}=\mE^\e\left|\mE^\e[\gamma^0_t|\mathscr{F}_t^{Y^\e}]\right|^{-2p}\leq\mE^\e\left[\mE^\e[|\gamma^0_t|^{-2p}|\mathscr{F}_t^{Y^\e}]\right]=\mE^\e[|\gamma^0_t|^{-2p}].
$$
Thus, the definition of $\gamma^0_t$ allows  us to obtain that
\ce
\mE^\e[|\gamma^0_t|^{-2p}]&=&\mE^\e\left[\exp\left\{-2p\int_0^t h(X_s^0)\dif
Y^{\e}_s+\frac{2p}{2}\int_0^t |h(X_s^0)|^2\dif s\right\}\right]\\
&=&\mE^\e\Bigg[\exp\left\{-2p\int_0^t h(X_s^0)\dif
Y^{\e}_s-\frac{4p^2}{2}\int_0^t |h(X_s^0)|^2\dif s\right\}\\
&&\bullet\exp\left\{\left(\frac{4p^2}{2}+\frac{2p}{2}\right)\int_0^t |h(X_s^0)|^2\dif s\right\}\Bigg]\\
&\leq&\exp\left\{(2p^2+p)CT\right\}\mE^\e\Bigg[\exp\left\{-2p\int_0^t h(X_s^0)\dif
Y^{\e}_s-\frac{4p^2}{2}\int_0^t |h(X_s^0)|^2\dif s\right\}\Bigg]\\
&=&\exp\left\{(2p^2+p)CT\right\},
\de
where the last step is based on the fact that  $\exp\left\{-2p\int_0^t h(X_s^0)\dif
Y^{\e}_s-\frac{4p^2}{2}\int_0^t |h(X_s^0)|^2\dif s\right\}$ is an exponential martingale under $\mP^\e$.

Similarly, we know that $\mE^\e(\gamma^\e_T)^2\leq\exp\left\{CT\right\}$. So, by
simple calculation, it holds that
$\mE\left|\rho^0_t(1)\right|^{-p}\leq\exp\left\{(2p^2+p+1)CT/2\right\}$. The proof is
complete.
\end{proof}

\bt\label{filcon} 
Suppose that {\bf Assume 1.-3.} hold.  Then it holds that for $\phi\in \cC_b^1(\mR^n)$
\be
\lim\limits_{\e\rightarrow0}\mE|\pi^{\e}_t(\phi)- \pi_t^0(\phi)|=0,
\label{filcones}
\ee
where $\cC_b^1(\mR^n)$ denotes the collection of all the functions which themselves and their first order partial derivatives are bounded and Borel measurable.
\et
\begin{proof}
For $\phi\in \cC^1_b(\mR^n)$, it follows from the H\"older inequality and Lemma \ref{es} that
\ce
\mE|\pi^{\e}_t(\phi)-
\pi_t^\e(\phi)|&=&\mE\left|\frac{\rho^{\e}_t(\phi)-\rho^{0}_t(\phi)}{\rho^0_t(1)}-\pi^{\e}_t(\phi)\frac{\rho^{\e}_t(1)-\rho^0_t(1)}{\rho^{0}_t(1)}\right|\\
&\leq&\mE\left|\frac{\rho^{\e}_t(\phi)-\rho^{0}_t(\phi)}{\rho^{0}_t(1)}\right|+\mE\left|\pi^{\e}_t(\phi)\frac{\rho^{\e}_t(1)-\rho^0_t(1)}{\rho^{0}_t(1)}\right|\\
&\leq&\left(\mE\left|\rho^{\e}_t(\phi)-\rho^{0}_t(\phi)\right|^{r_1}\right)^{1/r_1}\left(\mE\left|\rho^{0}_t(1)\right|^{-r_2}\right)^{1/r_2}\\
&&+\|\phi\|_{\cC^1_b(\mR^n)}\left(\mE\left|\rho^{\e}_t(1)-\rho^{0}_t(1)\right|^{r_1}\right)^{1/r_1}\left(\mE\left|\rho^{0}_t(1)\right|^{-r_2}\right)^{1/r_2}\\
&\leq&C\left(\mE\left|\rho^{\e}_t(\phi)-\rho^{0}_t(\phi)\right|^{r_1}\right)^{1/r_1}+C\|\phi\|_{\cC^1_b(\mR^n)}\left(\mE\left|\rho^{\e}_t(1)-\rho^{0}_t(1)\right|^{r_1}\right)^{1/r_1},
\de
where $1<r_1<2, r_2>1$ and $1/r_1+1/r_2=1$. 

Next, we estimate $\mE\left|\rho^{\e}_t(\phi)-\rho^{0}_t(\phi)\right|^{r_1}$. Note that
\ce
\mE\left|\rho^{\e}_t(\phi)-\rho^{0}_t(\phi)\right|^{r_1}&=&\mE^\e\left|\rho^{\e}_t(\phi)-\rho^{0}_t(\phi)\right|^{r_1}
\gamma^{\e}_T\leq (\mE^\e\left|\rho^{\e}_t(\phi)-\rho^{0}_t(\phi)\right|^{r_1p_1})^{1/p_1}
(\mE^\e(\gamma^\e_T)^{p_2})^{1/p_2}\\
&\leq&\exp\left\{CT\right\}(\mE^\e\left|\rho^{\e}_t(\phi)-\rho^{0}_t(\phi)\right|^{r_1p_1})^{1/p_1},
\de
where $1<p_1<2, 1<r_1p_1<2, p_2>1$ and $1/p_1+1/p_2=1$. And then we only need to observe $\mE^\e\left|\rho^{\e}_t(\phi)-\rho^{0}_t(\phi)\right|^{r_1p_1}$. Based
on the definitions of $\rho^{\e}_t(\phi), \rho^{0}_t(\phi)$ and the Jensen inequality, it holds that
\be
\mE^\e\left|\rho^{\e}_t(\phi)-\rho^{0}_t(\phi)\right|^{r_1p_1}&=&\mE^\e\left|\mE^\e[\phi(X_t^\e)\gamma^\e_t|\mathscr{F}_t^{Y^\e}]-\mE^\e[\phi(X^0_t)\gamma^0_t|\mathscr{F}_t^{Y^\e}]\right|^{r_1p_1}\no\\
&=&\mE^\e\left|\mE^\e[\phi(X_t^\e)\gamma^\e_t-\phi(X^0_t)\gamma^0_t|\mathscr{F}_t^{Y^\e}]\right|^{r_1p_1}\no\\
&\leq&\mE^\e\left[\mE^\e\left[\left|\phi(X_t^\e)\gamma^\e_t-\phi(X^0_t)\gamma^0_t\right|^{r_1p_1}\bigg|\mathscr{F}_t^{Y^\e}\right]\right]\no\\
&=&\mE^\e\left[\left|\phi(X_t^\e)\gamma^\e_t-\phi(X^0_t)\gamma^0_t\right|^{r_1p_1}\right]\no\\
&\leq&2^{r_1p_1-1}\mE^\e\left[\left|\phi(X_t^\e)\gamma^\e_t-\phi(X^0_t)\gamma^\e_t\right|^{r_1p_1}\right]\no\\
&&+2^{r_1p_1-1}\mE^\e\left[\left|\phi(X_t^0)\gamma^\e_t-\phi(X^0_t)\gamma^0_t\right|^{r_1p_1}\right]\no\\
&=:&I_1+I_2.
\label{i1i2}
\ee

First, we deal with $I_1$. By the H\"older inequality, it holds that
\be
I_1&\leq&\displaystyle
2^{r_1p_1-1}(\mE^\e\left[\left|\phi(X_t^\e)-\phi(X^0_t)\right|^{r_1p_1q_1}\right])^{1/q_1}
(\mE^\e\left|\gamma^\e_t\right|^{r_1p_1q_2})^{1/q_2} \no\\
&\leq&2^{r_1p_1-1}\|\phi\|_{\cC^1_b(\mR^n)}^{r_1p_1}(\mE^\e\left|X_t^\e-X^0_t\right|^{r_1p_1q_1})^{1/q_1}\no\\
&&\bullet \Bigg(\mE^\e\exp\left\{r_1p_1q_2\int_0^t
h(X_s^\e)\dif Y^{\e}_s-\frac{(r_1p_1q_2)^2}{2}\int_0^t|h(X_s^\e)|^2\dif s\right\}\no\\
&&\cdot\exp\left\{\frac{(r_1p_1q_2)^2}{2}\int_0^t|h(X_s^\e)|^2\dif s-\frac{r_1p_1q_2}{2}\int_0^t|h(X_s^\e)|^2\dif s\right\}
\Bigg)^{1/q_2}\no\\
&\leq&2^{r_1p_1-1}\|\phi\|_{\cC^1_b(\mR^n)}^{r_1p_1}(\mE^\e\left|X_t^\e-X^0_t\right|^{r_1p_1q_1})^{1/q_1}e^{\frac{r_1p_1}{2}(r_1p_1q_2-1)CT},
\label{i11}
\ee
where $1<q_1<2, 1<r_1p_1q_1<2, q_2>1$ and $1/q_1+1/q_2=1$, and the last step is based on the fact that the process
$\exp\left\{r_1p_1q_2\int_0^t h(X_s^\e)\dif Y^{\e}_s-\frac{(r_1p_1q_2)^2}{2}\int_0^t|h(X_s^\e)|^2\dif
s\right\}$ is an exponential martingale under $\mP^\e$. Note that 
\be
\mE^\e\left|X_t^\e-X^0_t\right|^{r_1p_1q_1}&=&\mE\left|X_t^\e-X^0_t\right|^{r_1p_1q_1}(\gamma^\e_T)^{-1}\no\\
&\leq& (\mE\left|X_t^\e-X^0_t\right|^2)^{r_1p_1q_1/2}\left(\mE(\gamma^\e_T)^{-2/(2-r_1p_1q_1)}\right)^{(2-r_1p_1q_1)/2}\no\\
&\leq&CR(\e)^{r_1p_1q_1/2},
\label{metrxzex}
\ee
where $R(\e):=\(C\frac{\e}{\delta_{\e}}+C(\delta_{\e}+1)\delta_{\e}+C\frac{(\delta_{\e}+1)\delta^2_{\e}}{\e}\)e^{CT}$ and the last step is based on Theorem \ref{xzerx}. Thus, by inserting (\ref{metrxzex}) in (\ref{i11}), we have that 
\ce
I_1\leq C\|\phi\|_{\cC^1_b(\mR^n)}^{r_1p_1}R(\e)^{r_1p_1/2}.
\de
We choose $\delta_\e$ as that in Remark \ref{del}, and obtain that $\lim\limits_{\e\rightarrow0}R(\e)=0$ and 
\be
\lim\limits_{\e\rightarrow0}I_1=0.
\label{i1es}
\ee

Next, for $I_2$, we know that
\ce
I_2&\leq&2^{r_1p_1-1}\|\phi\|_{\cC^1_b(\mR^n)}^{r_1p_1}\mE^\e\left[\left|\gamma^\e_t-\gamma^0_t\right|^{r_1p_1}\right]=2^{r_1p_1-1}\|\phi\|_{\cC^1_b(\mR^n)}^{r_1p_1}\mE\left[\left|\gamma^\e_t-\gamma^0_t\right|^{r_1p_1}\right](\gamma^\e_T)^{-1}\\
&\leq&2^{r_1p_1-1}\|\phi\|_{\cC^1_b(\mR^n)}^{r_1p_1}\left(\mE\left|\gamma^\e_t-\gamma^0_t\right|^2\right)^{r_1p_1/2}\left(\mE(\gamma^\e_T)^{-2/(2-r_1p_1)}\right)^{(2-r_1p_1)/2}\\
&\leq&C\left(\mE\left|\gamma^\e_t-\gamma^0_t\right|^2\right)^{r_1p_1/2}.
\de
Note that $\gamma^\e_t, \gamma^0_t$ have the following expressions
\ce
&&\gamma^\e_t=\exp\bigg\{\int_0^th(X_s^{\e})^i\dif U^{i}_s+\frac{1}{2}\int_0^t
\left|h(X_s^{\e})\right|^2\dif s\bigg\},\\
&&\gamma^0_t=\exp\bigg\{\int_0^th(X_s^0)^i\dif U^{i}_s+\int_0^th(X_s^0)^ih(X_s^{\e})^i\dif s-\frac{1}{2}\int_0^t
\left|h(X_s^0)\right|^2\dif s\bigg\}.
\de
So, by Theorem \ref{xzerx} and simple calculation, it holds that 
\ce
\lim\limits_{\e\rightarrow0}|\gamma^\e_t-\gamma^0_t|=0
\de
Moreover, ($\mathbf{H}_{h}$) admits us to get that
\ce
&&|\gamma^\e_t|^2\leq \exp\bigg\{\int_0^t 2h(X_s^{\e})^i\dif U^{i}_s-\frac{1}{2}\int_0^t
\left|2h(X_s^{\e})\right|^2\dif s\bigg\}e^{CT},\\
&&|\gamma^0_t|^2\leq \exp\bigg\{\int_0^t2h(X_s^0)^i\dif U^{i}_s-\frac{1}{2}\int_0^t
\left|2h(X_s^0)\right|^2\dif s\bigg\}e^{CT},
\de
and then $\mE|\gamma^\e_t|^2\leq e^{CT}, \mE|\gamma^0_t|^2\leq e^{CT}$. Thus, the dominated convergence theorem yields that 
\be
\lim\limits_{\e\rightarrow0}I_2=0.
\label{i2es}
\ee

Finally,  combining (\ref{i1es}) (\ref{i2es}) with (\ref{i1i2}), we obtain that 
\ce
\lim\limits_{\e\rightarrow0}\mE^\e\left|\rho^{\e}_t(\phi)-\rho^{0}_t(\phi)\right|^{r_1p_1}=0,
\de
and furthermore
$$
\lim\limits_{\e\rightarrow0}\mE\left|\rho^{\e}_t(\phi)-\rho^{0}_t(\phi)\right|^{r_1}=0.
$$
We recall that 
$$
\mE|\pi^{\e}_t(\phi)-
\pi_t^\e(\phi)|\leq C\left(\mE\left|\rho^{\e}_t(\phi)-\rho^{0}_t(\phi)\right|^{r_1}\right)^{1/r_1}+C\|\phi\|_{\cC^1_b(\mR^n)}\left(\mE\left|\rho^{\e}_t(1)-\rho^{0}_t(1)\right|^{r_1}\right)^{1/r_1}.
$$
Thus, taking the limit on two sides as $\e\rightarrow0$, one can get (\ref{filcones}). The proof is complete.
\end{proof}

\br
Here we can not give out the convergence rate of $\pi_t^\e$ to $\pi_t^0$. That is  because the convergence of the slow part to the homogenized system is in $L^2$ sense and is not in $L^p$ sense for any $p>1$.
\er

\section{Convergence of nonlinear filterings with correlated noises}\label{cornoi}

In the section, we study the nonlinear filtering problem of the system (\ref{Eq01}). First of all, we give out our assumption.

\begin{enumerate}[\bf{Assumption 4.}]
\item
\end{enumerate}
\begin{enumerate}[\bf{(i)}]
\item $\check{b}_1,\check{\sigma}_0, \check{\sigma}_1, \check{f}_1$ satisfy ($\mathbf{H}^1_{b_1, \sigma_1, f_1}$)-($\mathbf{H}^2_{b_1, \sigma_1, f_1}$), where $\check{b}_1, (\check{\sigma}_0, \check{\sigma}_1), \check{f}_1$ replace $b_1, \sigma_1, f_1$;
\end{enumerate}
\begin{enumerate}[\bf{(ii)}]
\item $\check{b}_2, \check{\sigma}_2, \check{f}_2$ satisfy ($\mathbf{H}^1_{b_2}$), ($\mathbf{H}^1_{\sigma_2}$) and ($\mathbf{H}^1_{f_2}$), respectively;
\end{enumerate}
\begin{enumerate}[\bf{(iii)}]
\item $\check{b}_2, \check{\sigma}_2, \check{f}_2$ satisfy ($\mathbf{H}^2_{\sigma_2}$), ($\mathbf{H}^1_{b_2, \sigma_2, f_2}$)-($\mathbf{H}^2_{b_2, \sigma_2, f_2}$), where $\check{b}_2, \check{\sigma}_2, \check{f}_2$ replace $b_2, \sigma_2, f_2$.
\end{enumerate}

Under {\bf Assumption 4.} {\bf (i)-(ii)}, by Theorem 1.2 in \cite{q2}, the system (\ref{Eq01}) has a unique strong solution denoted by $(\check{X}^\e_t,\check{Z}^\e_t)$. And then take any $x\in\mR^n$ and fix it. And consider the following SDE in $\mR^m$:
\ce\left\{\begin{array}{l}
\dif \check{Z}^x_t=\check{b}_2(x,\check{Z}^x_t)\dif t+\check{\sigma}_2(x,\check{Z}^x_t)\dif W_t+\int_{\mU_2}\check{f}_2(x,\check{Z}^x_t,u)\tilde{N}_{p_2}(\dif t, \dif u),\\
\check{Z}^x_0=\check{z}_0, \qquad t\geq0.
\end{array}
\right.
\de
Based on {\bf Assumption 4.} {\bf (ii)-(iii)}, it holds that the above equation has a unique invariant probability measure denoted as $\bar{\check{p}}(x,\dif z)$. So, set
\ce
\bar{\check{b}}_1(x):=\int_{\mR^m}\check{b}_1(x,z)\bar{\check{p}}(x,\dif z), 
\de
and by \cite[Lemma 3.1]{q00}, we know that $\bar{\check{b}}_1$ is Lipschitz continuous. So, we construct a SDE on the probability space $(\Omega, \mathscr{F}, \{\mathscr{F}_t\}_{t\in[0,T]}, \mP)$ as follows:
\be\left\{\begin{array}{l}
\dif \check{X}^0_t=\bar{\check{b}} _1(\check{X}^0_t)\dif t+\check{\sigma}_0(\check{X}^0_t)\dif B_t+\check{\sigma}_1(\check{X}^0_t)\dif V_t+\int_{\mU_1}\check{f}_1(\check{X}^0_{t-}, u)\tilde{N}_{p_1}(\dif t, \dif u),\\
\check{X}^0_0=\check{x}_0, \qquad\qquad 0\leq t\leq T.
\label{appequ2}
\end{array}
\right.
\ee
The solution of Eq.(\ref{appequ2}) is denoted as $\check{X}^0_t$. By the same deduction to that in Theorem \ref{xzerx}, we can obtain the following theorem.

\bt\label{corcon}
There exists a constant $C\geq 0$ independent of $\e, \delta_{\e}$ such that 
\ce
\mE\(\sup\limits_{0\leq t\leq T}|\check{X}^\e_t-\check{X}^0_t|^2\)\leq \(C\frac{\e}{\delta_{\e}}+C(\delta_{\e}+1)\delta_{\e}+C(\delta_{\e}+1)\frac{\delta^2_{\e}}{\e}\)e^{CT}.
\de
\et

\subsection{Nonlinear filtering problems with the system (\ref{Eq01})}

Next, for the observation process $\check{Y}^{\e}$ defined in (\ref{Eq01}), i.e.
\ce
\check{Y}_t^{\e}=V_t+\int_0^t\check{h}(\check{X}_s^{\e})\dif s+\int_0^t\int_{\mU_3}\check{f}_3(s,u)\tilde{N}_{\lambda}(\dif s, \dif u)+\int_0^t\int_{\mU\setminus\mU_3}\check{g}_3(s,u)N_{\lambda}(\dif s, \dif u),
\de
we assume:
\begin{enumerate}[\bf{Assumption 5.}]
\item
\end{enumerate}
\begin{enumerate}[\bf{(i)}]
\item $\check{h}$ is bounded and
$$
\int_0^T\int_{\mU_3}|\check{f}_3(s,u)|^2\nu_3(\dif u)\dif s<\infty.
$$
\end{enumerate}
\begin{enumerate}[\bf{(ii)}]
\item There exists a positive function $\check{L}(u)$ satisfying
\ce
\int_{\mU_3}\frac{\left(1-\check{L}(u)\right)^2}{\check{L}(u)}\nu_3(\dif u)<\infty
\de
such that $0<\check{l}\leq \check{L}(u)<\lambda(t,x,u)<1$ for $u\in\mU_3$, where $\check{l}$ is a constant.
\end{enumerate}

Now, denote
\ce
(\lambda^\e_t)^{-1}:&=&\exp\bigg\{-\int_0^t \check{h}^i(\check{X}_s^\e)\dif V^i_s-\frac{1}{2}\int_0^t
\left|\check{h}(\check{X}_s^\e)\right|^2\dif s-\int_0^t\int_{\mU_3}\log\lambda(s,\check{X}^{\e}_{s-},u)N_{\lambda}(\dif s, \dif u)\\
&&\quad\qquad -\int_0^t\int_{\mU_3}(1-\lambda(s,\check{X}^\e_s,u))\nu_3(\dif u)\dif s\bigg\}.
\de
Thus, by {\bf Assumption 5.} we know that $(\lambda^\e_t)^{-1}$ is an exponential martingale. Define a measure $\check{\mP}^\e$ via
$$
\frac{\dif \check{\mP}^\e}{\dif \mP}=(\lambda^\e_T)^{-1}.
$$
Under the probability measure $\check{\mP}^\e$, it follows from the Girsanov theorem that $\check{V}_t:=V_t+\int_0^t\check{h}(\check{X}_s^{\e})\dif s$ is a Brownian motion and $N_{\lambda}(\dif t,\dif u)$
is a Poisson random measure with the predictable compensator $\dif t\nu_3(\dif u)$. Moreover, by the same deduction to that in \cite[Lemma 3.1]{qd}, we know that $\lambda^\e_t$ satisfies the following equation
\be
\lambda^\e_t=1+\int_0^t\lambda^\e_s\check{h}(\check{X}^\e_s)^i\dif \check{V}^i_s+\int_0^t\int_{\mU_3}\lambda^\e_{s-}(\lambda(s,\check{X}^\e_{s-},u)-1)\tilde{N}(\dif s, \dif u),
\label{lme}
\ee
where $\tilde{N}(\dif s, \dif u):=N_{\lambda}(\dif t,\dif u)-\dif t\nu_3(\dif u)$. Set
\ce
&&\check{\rho}^\e_t(\psi):=\mE^{\check{\mP}^\e}[\psi(\check{X}^\e_t)\lambda^\e_t|\mathscr{F}_t^{\check{Y}^{\e}}],\\
&&\check{\pi}^{\e}_t(\psi):=\mE[\psi(\check{X}^{\e}_t)|\mathscr{F}_t^{\check{Y}^{\e}}], \qquad \psi\in\cB(\mR^n),
\de
where $\mE^{\check{\mP}^\e}$ stands for the expectation under the probability measure $\check{\mP}^\e$. And then by the Kallianpur-Striebel formula it holds that
\ce
\check{\pi}^{\e}_t(\psi)=\frac{\check{\rho}^\e_t(\psi)}{\check{\rho}^\e_t(1)}.
\de
In addition, we have the following result.

\bt(The Zakai equation)\label{zakequ}
For $\psi\in \cC^2_b(\mR^n)$, $\check{\rho}^{\e}_t(\psi)$ satisfies the following Zakai equation
\be
\check{\rho}^{\e}_t(\psi)&=&\check{\rho}^{\e}_0(\psi)+\int_0^t\check{\rho}^{\e}_s\(\big(\cL^{\check{X}^\e}\psi\big)(\cdot,\check{Z}^\e_s)\)\dif s+\int_0^t\left(\check{\rho}^{\e}_s(\psi \check{h}^i)+\check{\rho}^{\e}_s((\partial_j\psi)\check{\sigma}^{ji}_1)\right)\dif \check{V}^i_s\no\\
&&+\int_0^t\int_{\mU_3}\check{\rho}^{\e}_s\(\psi(\lambda(s,\cdot,u)-1)\)\tilde{N}(\dif s, \dif u),
\label{zakai1}
\ee
where
\ce
(\cL^{\check{X}^\e}\psi)(x,z)&:=&\frac{\partial \psi(x)}{\partial x_i}\check{b}^i_1(x,z)+\frac{1}{2}\frac{\partial^2\psi(x)}{\partial x_i\partial x_j}
(\check{\sigma}_0\check{\sigma}_0^T)^{ij}(x)+\frac{1}{2}\frac{\partial^2\psi(x)}{\partial x_i\partial x_j}
(\check{\sigma}_1\check{\sigma}_1^T)^{ij}(x)\\
&&+\int_{\mU_1}\Big[\psi\big(x+\check{f}_1(x,u)\big)-\psi(x)
-\frac{\partial \psi(x)}{\partial x_i}\check{f}^i_1(x,u)\Big]\nu_1(\dif u).
\de
\et
\begin{proof}
Applying the It\^o formula to $\psi(\check{X}^\e_t)$, one can have that
\ce
\psi(\check{X}^\e_t)&=&\psi(\check{x}_0)+\int_0^t(\cL^{\check{X}^\e}\psi)(\check{X}^\e_s,\check{Z}^\e_s)\dif s+\int_0^t(\nabla\psi)(\check{X}^\e_s)\check{\sigma}_0(\check{X}^\e_s)\dif B_s\\
&&+\int_0^t(\nabla\psi)(\check{X}^\e_s)\check{\sigma}_1(\check{X}^\e_s)\dif \left(\check{V}_s-\int_0^s\check{h}(\check{X}_r^{\e})\dif r\right)\\
&&+\int_0^t\int_{\mU_1}[\psi(\check{X}^\e_{s-}+\check{f}_1(\check{X}^\e_{s-}, u))-\psi(\check{X}^\e_{s-})]\tilde{N}_{p_1}(\dif s, \dif u).
\de
Note that $\lambda^\e_t$ satisfies (\ref{lme}). So, it follows from the It\^o formula that
\ce
\psi(\check{X}^\e_t)\lambda^\e_t&=&\psi(\check{x}_0)+\int_0^t\psi(\check{X}^\e_s)\lambda^\e_s\check{h}(\check{X}_s^{\e})^i\dif \check{V}^i_s\\
&&+\int_0^t\int_{\mU_3}\psi(\check{X}^\e_{s-})\lambda^\e_{s-}(\lambda(s,\check{X}^\e_{s-},u)-1)\tilde{N}(\dif s, \dif u)\\
&&+\int_0^t\lambda^\e_s(\cL^{\check{X}^\e}\psi)(\check{X}^\e_s,\check{Z}^\e_s)\dif s+\int_0^t\lambda^\e_s(\nabla\psi)(\check{X}^\e_s)\check{\sigma}_0(\check{X}^\e_s)\dif B_s\\
&&+\int_0^t\lambda^\e_s(\nabla\psi)(\check{X}^\e_s)\check{\sigma}_1(\check{X}^\e_s)\dif \left(\check{V}_s-\int_0^s\check{h}(\check{X}_r^{\e})\dif r\right)\\
&&+\int_0^t\int_{\mU_1}\lambda^\e_{s-}[\psi(\check{X}^\e_{s-}+\check{f}_1(\check{X}^\e_{s-}, u))-\psi(\check{X}^\e_{s-})]\tilde{N}_{p_1}(\dif s, \dif u)\\
&&+\int_0^t(\partial_j\psi)(\check{X}^\e_s)\check{\sigma}^{ji}_1(\check{X}^\e_s)\lambda^\e_s\check{h}(\check{X}_s^{\e})^i\dif s\\
&=&\psi(\check{x}_0)+\int_0^t\left(\psi(\check{X}^\e_s)\lambda^\e_s\check{h}(\check{X}_s^{\e})^i+\lambda^\e_s(\partial_j\psi)(\check{X}^\e_s)\check{\sigma}^{ji}_1(\check{X}^\e_s)\right)\dif \check{V}^i_s\\
&&+\int_0^t\int_{\mU_3}\psi(\check{X}^\e_{s-})\lambda^\e_{s-}(\lambda(s,\check{X}^\e_{s-},u)-1)\tilde{N}(\dif s, \dif u)\\
&&+\int_0^t\lambda^\e_s(\cL^{\check{X}^\e}\psi)(\check{X}^\e_s,\check{Z}^\e_s)\dif s+\int_0^t\lambda^\e_s(\nabla\psi)(\check{X}^\e_s)\check{\sigma}_0(\check{X}^\e_s)\dif B_s\\
&&+\int_0^t\int_{\mU_1}\lambda^\e_{s-}[\psi(\check{X}^\e_{s-}+\check{f}_1(\check{X}^\e_{s-}, u))-\psi(\check{X}^\e_{s-})]\tilde{N}_{p_1}(\dif s, \dif u).
\de
Taking the conditional expectation with respect to $\mathscr{F}_t^{\check{Y}^{\e}}$ under $\check{\mP}^\e$ on two hand sides of the above
equality, one could obtain  that
\ce
\mE^{\check{\mP}^\e}[\psi(\check{X}^\e_t)\lambda^\e_t|\mathscr{F}_t^{\check{Y}^{\e}}]&=&\mE^{\check{\mP}^\e}[\psi(\check{x}_0)|\mathscr{F}_t^{\check{Y}^{\e}}]\\
&&+\int_0^t\mE^{\check{\mP}^\e}[\left(\psi(\check{X}^\e_s)\lambda^\e_s\check{h}(\check{X}_s^{\e})^i+\lambda^\e_s(\partial_j\psi)(\check{X}^\e_s)\check{\sigma}^{ji}_1(\check{X}^\e_s)\right)|\mathscr{F}_t^{\check{Y}^{\e}}]\dif \check{V}^i_s\\
&&+\int_0^t\int_{\mU_3}\mE^{\check{\mP}^\e}[\psi(\check{X}^\e_{s-})\lambda^\e_{s-}(\lambda(s,\check{X}^\e_{s-},u)-1)|\mathscr{F}_t^{\check{Y}^{\e}}]\tilde{N}(\dif s, \dif u)\\
&&+\int_0^t\mE^{\check{\mP}^\e}[\lambda^\e_s(\cL^{\check{X}^\e}\psi)(\check{X}^\e_s,\check{Z}^\e_s)|\mathscr{F}_t^{\check{Y}^{\e}}]\dif s,
\de
where \cite[Theorem 1.4.7]{blr} is used. That is, it holds that
\ce
\check{\rho}^{\e}_t(\psi)&=&\check{\rho}^{\e}_0(\psi)+\int_0^t\check{\rho}^{\e}_s\(\big(\cL^{\check{X}^\e}\psi\big)(\cdot,\check{Z}^\e_s)\)\dif s+\int_0^t\left(\check{\rho}^{\e}_s(\psi \check{h}^i)+\check{\rho}^{\e}_s((\partial_j\psi)\check{\sigma}^{ji}_1)\right)\dif \check{V}^i_s\no\\
&&+\int_0^t\int_{\mU_3}\check{\rho}^{\e}_s\(\psi(\lambda(s,\cdot,u)-1)\)\tilde{N}(\dif s, \dif u).
\de
The proof is complete.
\end{proof}

In the following, we define the nonlinear filtering of $\check{X}_t^0$ with respect to $\mathscr{F}_t^{\check{Y}^{\e}}$. Set
\ce
\lambda^0_t&:=&\exp\bigg\{\int_0^t\check{h}^i(\check{X}_s^0)\dif \check{V}^{i}_s-\frac{1}{2}\int_0^t
\left|\check{h}(\check{X}_s^0)\right|^2\dif s+\int_0^t\int_{\mU_3}\log\lambda(s,\check{X}^0_{s-},u)N_{\lambda}(\dif s, \dif u)\\
&&\quad\qquad +\int_0^t\int_{\mU_3}(1-\lambda(s,\check{X}^0_s,u))\nu_3(\dif u)\dif s\bigg\},
\de
and furthermore
\ce
\check{\rho}^0_t(\psi)&:=&\mE^{\check{\mP}^\e}[\psi(\check{X}^0_t)\lambda^0_t|\mathscr{F}_t^{\check{Y}^{\e}}],\\
\check{\pi}^0_t(\psi)&:=&\frac{\check{\rho}^0_t(\psi)}{\check{\rho}^0_t(1)}.
\de
And then by the similar deduction to that in Theorem \ref{zakequ}, it holds that $\check{\rho}^0_t$ satisfies the following equation
\be
\check{\rho}^0_t(\psi)&=&\check{\rho}^0_0(\psi)+\int_0^t\check{\rho}^0_s\(\cL^{\check{X}^0}\psi\)\dif s+\int_0^t\check{\rho}^0_s\Big(((\partial_j\psi)\check{\sigma}^{ji}_1)(\check{h}^i(\cdot)-\check{h}^i(\check{X}_s^\e))\Big)\dif s\no\\
&&+\int_0^t\left(\check{\rho}^0_s(\psi \check{h}^i)+\check{\rho}^0_s((\partial_j\psi)\check{\sigma}^{ji}_1)\right)\dif \check{V}^i_s\no\\
&&+\int_0^t\int_{\mU_3}\check{\rho}^0_s\(\psi(\lambda(s,\cdot,u)-1)\)\tilde{N}(\dif s, \dif u),
\label{zakai2}
\ee
where 
\ce
(\cL^{\check{X}^0}\psi)(x)&:=&\frac{\partial \psi(x)}{\partial x_i}\bar{\check{b}}^i_1(x)+\frac{1}{2}\frac{\partial^2\psi(x)}{\partial x_i\partial x_j}
(\check{\sigma}_0\check{\sigma}_0^T)^{ij}(x)+\frac{1}{2}\frac{\partial^2\psi(x)}{\partial x_i\partial x_j}
(\check{\sigma}_1\check{\sigma}_1^T)^{ij}(x)\\
&&+\int_{\mU_1}\Big[\psi\big(x+\check{f}_1(x,u)\big)-\psi(x)
-\frac{\partial \psi(x)}{\partial x_i}\check{f}^i_1(x,u)\Big]\nu_1(\dif u).
\de

\subsection{The relationship of $\check{\pi}^\e$ and $\check{\pi}^0$} 

\subsubsection{The case of $\check{f}_3=\check{g}_3=0$}

In the case of $\check{f}_3=\check{g}_3=0$, by the similar deduction to that in Theorem \ref{filcon}, one can obtain the following result.

\bt
Suppose that {\bf Assume 4.-5.} hold.  Then it holds that for $\phi\in \cC_b^1(\mR^n)$,
\ce
\lim\limits_{\e\rightarrow0}\mE|\check{\pi}^{\e}_t(\phi)- \check{\pi}_t^0(\phi)|=0.
\de
\et

\subsubsection{The case of $\check{f}_3\neq\check{g}_3\neq0$}

In the case of $\check{f}_3\neq\check{g}_3\neq0$, we first prepare two important lemmas. Since their proofs are similar to that of \cite[Lemma 5.1, 5.2]{q00}, we omit them.

\bl\label{rho1}
Under {\bf Assumption 4.-5.}, it holds that for any $t\in[0,T]$,
\be
(\check{\rho}^0_t(1))^{-1}<\infty, \qquad \mP\ a.s..
\label{roes}
\ee
\el

\bl\label{tigh}
Under {\bf Assumption 4.-5.}, $\{\check{\rho}^{\e}_t, t\in[0,T]\}$ is relatively weakly compact in $D([0,T], \cM(\mR^n))$, where $\cM(\mR^n)$ denotes the set of bounded Borel measures on $\mR^n$.
\el

In the following, we assume more:

\begin{enumerate}[\bf{Assumption 6.}]
\item
\end{enumerate}
\begin{enumerate}[]  
\item $\{\check{Z}_{\e t}^{\e}, t\in[0,T]\}$ is tight.
\end{enumerate}

Now, we state and prove the main theorem in the section.

\bt\label{filcon2}
Suppose that {\bf Assumption 4.-6.} hold. Then $\check{\pi}^{\e}_t$ converges weakly to $\check{\pi}^0_t$ as $\e\rightarrow0$ for any $t\in[0,T]$.
\et
\begin{proof}
By the definition of $\check{\pi}^{\e}_t$, $\check{\pi}^0_t$, it holds that for $\phi\in\cC^2_b(\mR^n)$, 
\ce
\check{\pi}^{\e}_t(\phi)-\check{\pi}^0_t(\phi)=\frac{\check{\rho}^{\e}_t(\phi)-\check{\rho}^0_t(\phi)}{\check{\rho}^0_t(1)}-\check{\pi}^{\e}_t(\phi)\frac{\check{\rho}^{\e}_t(1)-\check{\rho}^0_t(1)}{\check{\rho}^0_t(1)}.
\de
So, in order to prove $\check{\pi}^{\e}_t(\phi)-\check{\pi}^0_t(\phi)$ converges weakly to $0$, in terms of Lemma \ref{rho1}, we only need to show that $\check{\rho}^{\e}_t(\phi)-\check{\rho}^0_t(\phi)$ converges weakly to $0$ as $\e\rightarrow0$. Besides, it follows from Lemma \ref{tigh} that there exist a weak convergence subsequence $\{\check{\rho}^{\e_k}_t, k\in\mN\}$ and a measure-valued process $\bar{\check{\rho}}_t$ such that $\check{\rho}^{\e_k}_t$ converges weakly to $\bar{\check{\rho}}_t$ as $k\rightarrow\infty$. Therefore, we just need to prove that for $t\in[0,T]$, $\bar{\check{\rho}}_t(\phi)-\check{\rho}^0_t(\phi)$ converges weakly to $0$ as $\e\rightarrow0$.

Next, we search for the equations which $\bar{\check{\rho}}_t(\phi)$ solves. By Theorem \ref{corcon} and (\ref{zakai1}), we follow up the line of \cite[Theorem 5.3]{q00} and obtain that $\bar{\check{\rho}}_t(\phi)$ 
satisfies the following equation
\be
\bar{\check{\rho}}_t(\phi)&=&\bar{\check{\rho}}_0(\phi)+\int_0^t\bar{\check{\rho}}_s\(\cL^{\check{X}^0}\phi\)\dif s+\int_0^t\left(\bar{\check{\rho}}_s(\phi\check{h}^i)+\bar{\check{\rho}}_s((\partial_j\phi)\check{\sigma}^{ji}_1)\right)\dif \check{V}^i_s\no\\
&&+\int_0^t\int_{\mU_3}\bar{\check{\rho}}_s\(\phi(\lambda(s,\cdot,u)-1)\)\tilde{N}(\dif s, \dif u).
\label{zakai3}
\ee

Besides, by (\ref{zakai2}) and Theorem \ref{corcon}, we know that, there exists a  measure-valued process $\bar{\check{\rho}}^0_t$ such that $\check{\rho}^0_t(\phi)$ converges to $\bar{\check{\rho}}^0_t$ $\mP$ a.s. and $\bar{\check{\rho}}^0_t$ satisfies the following equation
\be
\bar{\check{\rho}}^0_t(\phi)&=&\bar{\check{\rho}}^0_0(\phi)+\int_0^t\bar{\check{\rho}}^0_s\(\cL^{\check{X}^0}\phi\)\dif s+\int_0^t\left(\bar{\check{\rho}}^0_s(\phi \check{h}^i)+\bar{\check{\rho}}^0_s((\partial_j\phi)\check{\sigma}^{ji}_1)\right)\dif \check{V}^i_s\no\\
&&+\int_0^t\int_{\mU_3}\bar{\check{\rho}}^0_s\(\phi(\lambda(s,\cdot,u)-1)\)\tilde{N}(\dif s, \dif u).
\label{zakai4}
\ee

Note that Eq.(\ref{zakai3}) and Eq.(\ref{zakai4}) are the same. Thus, it follows from \cite[Theorem 3.9]{q3} that for any $t\in[0,T]$,
\ce
\bar{\check{\rho}}_t=\bar{\check{\rho}}^0_t, \quad \mP.a.s.
\de
That is, $\bar{\check{\rho}}_t(\phi)-\check{\rho}^0_t(\phi)$ converges weakly to $0$ as $\e\rightarrow0$. The proof is complete.
\end{proof}

\section{Conclusion}\label{con}

In the paper, we consider nonlinear filtering problems of multiscale systems in two cases-correlated sensor L\'evy noises and correlated L\'evy noises. First of all, we prove that the slow part of the origin system converges to the homogenized system in the uniform mean square sense. Next, in the case of correlated sensor L\'evy noises, the nonlinear filtering of the slow part is shown to approximate that of the homogenized system in $L^1$ sense. However, in the case of correlated L\'evy noises, we prove that the nonlinear filtering of the slow part converges weakly to that of the homogenized system.

\bigskip

\textbf{Acknowledgements:}

The author would like to thank Professor Xicheng Zhang for his valuable discussions. The author also thanks Professor Renming Song for providing her an excellent environment to work in the University of Illinois at Urbana-Champaign.

\end{document}